\documentclass[12pt,twoside,reqno]{amsart}
\usepackage[colorlinks=true,citecolor=blue]{hyperref}
\usepackage{mathptmx, amsmath, amssymb, amsfonts, amsthm, mathptmx, enumerate, color}
\usepackage{graphicx}
\usepackage{epstopdf}

\usepackage[ngerman,english]{babel}
\usepackage[utf8]{inputenc}
\usepackage{bbm}
\usepackage{subcaption}
\usepackage{algorithm2e}

\DeclareMathOperator{\diag}{diag}

\DeclareMathOperator{\rank}{rank}

\theoremstyle{definition}

\numberwithin{equation}{section}

\begin{document}
\setcounter{page}{1}

\vspace*{1.0cm}
\title[Immunity to Increasing Condition Numbers of Linear Superiorization
versus Linear Programming]
{\textbf{Immunity to Increasing Condition Numbers of Linear Superiorization
		versus Linear Programming}}
\author[J. Schröder, Y. Censor, P. Süss and K.-H.
Küfer]{ Jan Schröder$^{1}$, Yair Censor$^{2}$, Philipp Süss$^{1,*}$ and Karl-Heinz
	Küfer$^{1}$}
\maketitle
\vspace*{-0.6cm}

\begin{center}
{\footnotesize {\it

$^{1}$Optimization Department, Fraunhofer - ITWM, Kaiserslautern
67663, Germany\\
$^{2}$Department of Mathematics, University of Haifa, Mt. Carmel,
Haifa 3498838, Israel\\
}}

\vskip 4mm
\small{March 5, 2026. Revised: June 19, 2026}

\end{center}

\vskip 4mm {\small\noindent {\bf Abstract.}
Given a family of linear constraints and a linear objective function
one can consider whether to apply a Linear Programming (LP) algorithm
or use a Linear Superiorization (LinSup) algorithm on this data. In
the LP methodology one aims at finding an optimal point, i. e., a point
that fulfills the constraints and has the minimal value of the objective
function over these constraints.

The Linear Superiorization approach considers the same data as in linear
programming problems but instead of attempting to solve with
linear programming methods it employs perturbation resilient feasibility-seeking
algorithms that steer the iterations toward a feasible point with reduced (not
necessarily minimal) objective function value. This aim of the superiorization
method (SM) is less demanding than aiming to reach full-fledged constrained
optimality and it places more importance on reaching feasibility than
on reaching optimality.

Previous studies (e.g., \cite{linsup}) compared LP and LinSup in terms of their respective
outputs and the resources they use. Here, we investigate classical LP approaches and LinSup in terms of their sensitivity to condition numbers
of the system of linear constraints. 

Condition numbers are a measure
for the impact of deviations in the input data on the output of a
problem and, in particular, they describe the factor of error propagation
when given wrong or erroneous data. 

Therefore, the ability of LP and
LinSup to cope with increased condition numbers, thus with ill-posed
problems, is an important matter to consider which was not studied
until now. We investigate experimentally the advantages and disadvantages
of both LP and LinSup on exemplary sets of data of problems of linear programming
with multiple condition numbers and different problem dimensions.

 }

\renewcommand{\thefootnote}{}
\footnotetext{ $^*$Corresponding author.
\par
E-mail address: philipp.suess@itwm.fraunhofer.de (P. Süss).
}

\section{Introduction}\label{sec:Introduction}

In this paper we investigate experimentally and discuss the behavior
of the linear superiorization (LinSup) methodology for various different
condition numbers and compare it to the behavior of several linear
programming (LP) solvers. The paper is in the spirit of \cite{linsup} and contains a follow-up analysis regarding the impact of condition numbers on the runtime of the respective algorithms.

We set up linear programming problems with
different condition numbers of the linear constraints in a manner
that keeps the problems structurally similar in order to enable the
comparative experimental study. Then we apply to these problems the
LinSup superiorization method and some established implementations
of the simplex algorithm and an interior point method. 

Imposing stopping
rules that are fair toward all algorithms, we compare the outputs
on problems with increasing condition numbers and on problems with
increasing dimensions. 

\subsection{The experimental methodology}\label{subsec:The-experimental-methodology}

LinSup is not an LP solver. As an offspring of the superiorization
methodology (SM), LinSup perturbs linear feasibility-seeking algorithms
to reduce objective function values while retaining the feasibility-seeking
nature of the algorithm and without paying a high computational price.

The declared aim of LinSup is to find a feasible point whose objective
function value is ``superior'' (meaning smaller or equal) to that
of a feasible point that can be reached by the corresponding unperturbed
feasibility-seeking algorithm, i. e., the exactly same feasibility-seeking
algorithm that the LinSup employs. This aim is less demanding than
full-fledged linear programming but more demanding than plain feasibility-seeking.\footnote{Although we came across it after the SM was developed, one can borrow
	support for the reasoning of the SM from the American scientist and
	Nobel-laureate Herbert Simon who was in favor of “satisficing” rather
	than “maximizing”. Satisficing is a decision-making strategy that
	aims for a satisfactory or adequate result, rather than the optimal
	solution. This is because aiming for the optimal solution may necessitate
	needless expenditure of time, energy and resources. The term “satisfice”
	was coined by Herbert Simon in 1956 (\cite{simon}). See also: https://en.wikipedia.org/wiki/Satisficing.}

In view of this difference between LinSup and an LP solver, but recognizing
that both are initialized at a starting point which is not known to
be a feasible point, we maintain the following experimental methodology.
For a given problem we run LinSup until it reaches some predetermined
fixed small threshold level of ``infeasibility'' (i. e., ``constraints
violation''). When LinSup stopped according to this stopping rule
we record the objective function value at this point and the runtime
elapsed. 

Then we run on the same problem, starting from the same initialization
point, the LP solver that we wish to include in our comparative study,
and stop it when it reaches the same predetermined fixed level of
``infeasibility'' that was used for stopping the LinSup run. When
stopped we record the objective function value at this point and the
runtime elapsed for this LP solver as well.

We also record the infeasibility and the objective function values
along the whole runtime and use these in order to demonstrate the
behavior of each algorithm by plotting its objective function values
and its infeasibility values versus time, respectively. More details
are supplied later in the paper.

\subsection{Context and previous work}\label{subsec:Context-and-previous}

The superiorization methodology (SM) is a relatively new method that
aims to improve the performance of an existing iterative algorithm
by interlacing into it perturbation steps (\cite{censor2022superiorization}).
In the context of optimization problems, this concept can be applied
by interlacing into a feasibility-seeking algorithm (for example a
projection method) perturbations of negative gradient steps that reduce
the objective function values. 

Since its development, the SM has successfully
been applied in various practical applications such as intensity-modulated
radio therapy (\cite{supimrt}), image reconstruction (\cite{supimagerec})
and telecommunication networks (\cite{supfink}). Recent works include the use of machine learning for the determination of the perturbation steps (\cite {jia}) and systematic biasing of steepest descent methods in order to obtain solutions that are less dependent on noise in the data (\cite{nikazad}). Both have been applied to imaging.

In the particular context of linear programming problems (LPs), there is a vast literature
for both feasibility-seeking problems (for example the Agmon-Motzkin-Schönberg
(AMS) algorithm \cite{agmon,motzkinschoenberg},\cite[algorithm 5.4.1]{parallelopt})
and for the LP problem (e.g., simplex algorithm, ellipsoid method
\cite{lps}). 

Furthermore, for the LinSup case, the ``guarantee question of the
SM'' has been answered positively. This is the question whether superiorization
can actually converge to a feasible point with a smaller or equal
objective function value than that of the point reached by the un-superiorized
version of the same feasibility-seeking algorithm. In \cite{supconcmeasure}
the authors employ the principle of concentration of measure to show
that this result holds with high probability. A similar conclusion
for the nonlinear case is yet to be found.

In this paper we compare experimentally the SM with optimization
algorithms for LPs in an organized, reproducible and fair manner for
problems of varying difficulty, indicated by the problems' dimensions
and their condition numbers. 

Accordingly, our work is an extension of the results in \cite{linsup},
where other features of LinSup were compared to the simplex method. We augment these
results by varying the condition numbers across multiple problem instances.
The condition number of the linear system has been known to have a
significant impact on the performance of certain methods, often leading
to the failure of an algorithm (\cite{condnrmethods}). Since many
real-world problems have an inherently high condition number, immunity
against ill-conditioning of problems is a desirable property for any
algorithm. Often it is enough to have a rough understanding of the
order of magnitude of the condition number and there are efficient
algorithms for its estimation (\cite{condnrestimate},\cite{condtriest}).
For literature on preconditioning of matrix problems see, e.g., \cite{precondreview}.

\subsection{Our results}\label{subsec:Our-results}

In this project we used LP data aiming at the task of finding a feasible
point with reduced (not necessarily minimal) objective function value.
In our numerical work the superiorization method LinSup performed
better than LP solvers for this task. This better performance, for
this task, was investigated for problems with increasing condition
numbers of the linear systems and with increasing problem sizes. We
found that in each of these cases LinSup delivered better results
than some established LP solvers.

\subsection{Outline}\label{subsec:Outline}

We give a brief overview of the superiorization method and of condition
numbers in Sections \ref{superiorization} and \ref{condnr}, respectively.
We set up multiple linear programming problems with varying condition
numbers in Section \ref{formulation} and run the LinSup method and
several LP algorithms on them.

For details of our implementation of the LinSup method see Section
\ref{implDetails}. We give a brief description of the optimization
suites that we use for our comparisons in Section \ref{optalgs} and
of the experimental setup in Section \ref{expsetup}. In Section \ref{results}
we present the numerical results, followed by a discussion of remaining
challenges and future work in Section \ref{conclusion}.

\section{Preliminary background}\label{sec:Preliminaries}

\subsection{The superiorization methodology}

\label{superiorization}

The Superiorization Methodology (SM) evolved from the investigation
of feasibility-seeking models of some important real-world problems
such as image reconstruction from projections and radiation therapy
treatment planning. 

Feasibility-seeking algorithms, mainly projection
methods, generate iterative sequences that (under reasonable conditions)
converge to a point in the feasible set. Their main advantage is that
they perform projections onto the individual sets whose intersection
is the feasible set and not directly onto the feasible set and the
underlying premise is that such projections onto the individual sets
are more manageable.

When one wishes to find feasible points with a reduced, not necessarily
minimized, value of an imposed objective function then the SM comes
into play. The principle of the SM is to interlace into the iterates
of a feasibility-seeking iterative process perturbations that steer
the iterates toward superior (meaning smaller or equal) objective
function values without losing the overall convergence of the sequence
of perturbed iterates to a feasible point. To this end ``bounded
perturbations'' are used.

How all this is done is rigorously described in earlier papers on
the SM, consult, e.g., \cite{censor2022superiorization} for a recent
review, read also \cite{supprobstructs}. A key feature of the SM
is that it does not aim for a constrained optimal function value,
but is content with settling for a feasible point with reduced objective
function value -- reduced in comparison to the objective function
value of a feasible point that would be reached by the same feasibility-seeking
algorithm without perturbations. 

For many applications this is sufficient, in particular, whenever
the introduction of an objective function is only a secondary goal.
Fulfillment of the constraints, in this context, is considered by
the modeler of the real-world problem to be much more important, see,
e.g., \cite{censor2022superiorization,supimagerec,supimrt}.

Many papers on the SM are cited in \cite{supwebpage} which is a Webpage
dedicated to superiorization and perturbation resilience of algorithms
that contains a continuously updated bibliography on the subject.
This Webpage (which includes 208 entries as of June 6, 2026) is a
source for the wealth of work done in this field to date, including
two special issues of journals \cite{suptheoryapp} and \cite{supvsopt}
dedicated to research of the SM. Recent work includes \cite{bonackersup,abbasi,supvsreg,imagerecproj,supbregmanop,suptomography}. 

We find especially interesting the recent work of Ma et al. \cite{sahinidis}
who proposed a novel decomposition framework for derivative-free optimization
(DFO) algorithms which significantly extended the scope of current
DFO solvers to larger-scale problems. They proved that their proposed
framework closely relates to the superiorization methodology. 

Information
about feasibility-seeking algorithms and the convex feasibility problem,
that lie at the foundation of the SM approach, can be found, e.g.,
in Bauschke and Borwein's SIAM review paper \cite{bb96} and in \cite{swiss}.

\subsection{Condition numbers}

\label{condnr} The relative condition number is a measure of the
impact of deviations in the input data on the output data of a problem.
In particular, it describes the factor of error propagation when given
wrong or erroneous data. 

Let the function $f\colon\mathbb{R}^{n}\to\mathbb{R}^{m}$
represent some mathematical problem and let $x\in\mathbb{R}^{n}$
be its input, where $\mathbb{R}^{n}$ stands for the Euclidean $n$-dimensional
space. Denote with $\tilde{x}\in\mathbb{R}^{n}$ the disturbed input
data. Then the relative condition number of the problem at the point
$x$ is defined as (see \cite[equation (12.4)]{condnrtrefethen})
\begin{equation}
	\kappa_{\text{rel}}:=\limsup_{\tilde{x}\rightarrow x}\frac{\Vert f(\tilde{x})-f(x)\Vert}{\Vert\tilde{x}-x\Vert}\frac{\Vert x\Vert}{\Vert f(x)\Vert},
\end{equation}
as long as $f(x)\neq0$. Otherwise, it is $\kappa_{\text{rel}}=\infty$.

In particular, the condition number is independent of a chosen numerical
algorithm for solving the problem $f$, but the algorithms convergence
speed may depend on the magnitude of the condition number (see \cite[section 4]{condnrmethods}).

In the following we are interested in the condition number of matrices.
When $f(x)=Ax$ describes the problem of matrix multiplication, where
$A\in\mathbb{R}^{m\times n}$ and with $\Vert\cdot\Vert=\Vert\cdot\Vert_{2}$,
the above formula becomes

\begin{equation}
	\kappa_{\text{rel}}=\limsup_{\tilde{x}\rightarrow x}\frac{\Vert A(\tilde{x}-x)\Vert}{\Vert\tilde{x}-x\Vert}\frac{\Vert x\Vert}{\Vert Ax\Vert}.
\end{equation}

Since $f$ is differentiable and by writing $\tilde{x}-x=hv$, for
some unit vector $v$ and $h=\Vert\tilde{x}-x\Vert$ we get 
\begin{equation}
	\kappa_{\text{rel}}=\lim_{h\rightarrow0}\frac{\Vert A(x+hv-x)\Vert}{h}\frac{\Vert x\Vert}{\Vert Ax\Vert}=\Vert A\Vert\frac{\Vert x\Vert}{\Vert Ax\Vert}\leq\Vert A\Vert\Vert A^{-1}\Vert.
\end{equation}
The term on the right is called the ``condition number of the matrix
$A$'' (see \cite[Equation (12.15)]{condnrtrefethen}), where $A^{-1}$
denotes the inverse or, if $A$ is non-square, the pseudo-inverse
of $A$, 
\begin{equation}
	\kappa(A):=\Vert A\Vert\Vert A^{-1}\Vert=\frac{\sigma_{\text{max}}}{\sigma_{\text{min}}},
\end{equation}
where $\sigma_{\text{max}}$ and $\sigma_{\text{min}}$ are the maximal
and minimal nonzero singular values of $A$, respectively. 

The condition number
plays a significant role in the analysis of numerical problems and
is subject to extensive studies in the literature (\cite{condnrdeuflhard},
\cite{condnrtrefethen}, \cite{condnrestimate}, \cite{condnrmethods}).
Several methods exist to improve high condition numbers of ill-conditioned
problems (these are the, so-called, pre-conditioning methods, see
e.g., \cite{precondreview}) in order to increase the accuracy of
calculated solutions. This is often necessary because many real-world
applications give rise to condition numbers of significant magnitude.
This is the key motivation for the investigation in this paper.

\section{Problem Setup and Implementation Details}

\label{setup} 

\subsection{The problem formulation}

\label{formulation} We consider LP problems of the form

\begin{alignat}{2}
	& \!\min_{x\in\mathbb{R}^{n}} &  & \ \langle c,x\rangle\nonumber \\
	& \text{s.t. } &  & Ax\leq b\label{eq:LP}\\
	&  &  & \ell\leq x\leq u,\nonumber 
\end{alignat}
where $A\in\mathbb{R}^{m\times n},b\in\mathbb{R}^{m},c,\ell,u\in\mathbb{R}^{n}$.
Write $A=U\Sigma V$ via the singular value decomposition (cf. \cite{svd})
with semi-orthogonal matrices $U$ and $V$ and diagonal matrix $\Sigma=\diag(\sigma_{1},\dots,\sigma_{q})$
of singular values. Without loss of generality, let $\sigma_{1}\geq\ \dots\ \geq\sigma_{q}>0$.
The condition number $\kappa$ of $A$ is (see \cite{condnrdeuflhard})
\begin{equation}
	\kappa(A):=\frac{\sigma_{\text{max}}}{\sigma_{\text{min}}}=\frac{\sigma_{1}}{\sigma_{q}}.
\end{equation}
We want to construct a sequence of matrices $A_{\kappa}$ of a specified
condition number $\kappa$ in such a way that they remain structurally
similar to each other. To this end we reverse the singular value decomposition,
that is, we create exactly one pair of $U$ and $V$ which contains
the structure of the problem and construct, for different values of
$\kappa$, diagonal matrices $\Sigma_{\kappa}$, which impose the
condition number of the problem. Then, we calculate $A_{\kappa}:=U\Sigma_{\kappa}V$. 

As any matrix has a singular value decomposition, this makes it possible
to define any matrix via this approach too. In our construction we
focus, without loss of generality, on matrices of full rank, because
otherwise one can remove rows or columns until full rank is achieved.
For any chosen value of $\kappa$ we construct the diagonal matrix
$\Sigma_{\kappa}$, by setting $q:=\min(m,n)=\rank(A)$ and defining
singular values as
\begin{equation}
	\sigma_{i}:=\frac{t}{z_{i}}+\frac{1-t}{s},
\end{equation}\label{singularValueGenerator}
where $t:=\frac{\kappa-1}{q-1}$, $s=10$ and $z_{i}=\frac{s\cdot i}{q}$.
We chose this setup of the singular values $\sigma_{i}$, because
in real-world applications the singular values often behave approximately
proportional to $\frac{1}{i}$ (instead of linear which seems like
an obvious first choice for our problem). 

The parameter $s$ is used to control the magnitude of the singular
values since in this model we always have $\sigma_{q}=\frac{1}{s}$.
Due to the choice of $t$ it is easy to see that 
\begin{equation}
	\frac{\sigma_{1}}{\sigma_{q}}=\frac{\sigma_{\text{max}}}{\sigma_{\text{min}}}=\kappa
\end{equation}
as desired. For $U$ and $V$ we generate random semi-orthogonal matrices
and set 
\begin{equation}
	A=U\cdot\diag(\sigma_{1},\dots,\sigma_{q})\cdot V. 
\end{equation}

Furthermore, in order to guarantee the non-emptiness of the feasible
set in (\ref{eq:LP}) we set $b:=A\mathbbm{1}+\mathbbm{1}$, $u=-\ell=100\cdot\mathbbm{1}$
and randomly choose $c$. This choice of parameters implies the feasibility
of $x=\mathbbm{1}$. 

We use the above set up to generate LP problems for multiple dimensions,
with matrix sizes ranging from $80\times100$ to $4000\times5000$
(cf. \cite{linsup}).

\subsection{The superiorization algorithm}

\label{implDetails}

We aim to apply separately LP solvers and the LinSup method for the
data $A=(a^{i})_{i=1}^{m},b,c,\ell,u$ that appears in \eqref{eq:LP}.
In LinSup, presented here in Algorithm \ref{supalg}, we chose for
the feasibility-seeking algorithm the projection method of Agmon Motzkin
and Schönberg (AMS) as the ``basic
algorithm''. 

For some (usually small value) $r\geq0$, this algorithm sequentially cyclically projects onto
the individual half-spaces $\langle a^{i},x\rangle\leq b_{i}$ using
(see, e.g., \cite[p. 6]{linsup}) for each individual projection
the closed-form formula
\begin{equation}
	P_{i}(x):=\begin{cases}
		x-\lambda\frac{\langle a^{i},x\rangle-b_{i}+r\Vert a^{i}\Vert}{\Vert a^{i}\Vert^{2}}a^{i}\text{, if \ensuremath{\langle a^{i},x\rangle>b_{i},}}\\
		x\text{, otherwise,}
	\end{cases}\label{eq:amsstep}
\end{equation}
and a full sweep through all half-spaces is done by the algorithmic
operator $\mathcal{A}$, which is a composition of the individual
projections 
\begin{equation}
	\mathcal{A}:=P_{m}\circ\dots\circ P_{1}.\label{eq:amssweep}
\end{equation}
The parameter $\lambda$ describes a relative overshoot or undershoot of the projection. For our purposes we set $\lambda=1$ so that (if $r=0$) the projection $P_i$ maps precisely onto the $i$-th individual half-space $\langle a^i,x\rangle\leq b_i$. 

For the direction vectors in the perturbations used in LinSup of Algorithm
\ref{supalg} we chose the normalized negative gradient of the objective
function in \eqref{eq:LP}, which is constant throughout and equals
$-\frac{c}{\Vert c\Vert}$. 

For the step-sizes $\eta_{k}$ we take an exponentially decreasing
null sequence $(10\cdot0.99^{n})_{n\in\mathbb{N}}$ with restarts
every $\tau_{\text{restart}}=20$ iterations as described in \cite[p. 6]{restarts}.
The starting step-size $\eta_{0}=10$ is decreased by the kernel $\alpha=0.99$
in each iteration, unless there is a restart. In that case, we set
$\eta_{k}=\eta_{0}\alpha^{\rho}$, where $\rho$ is the number of
previous restarts during this run.

In other words, the $k$-th iteration consists of a gradient step
$-\eta_{k}\frac{c}{\Vert c\Vert}$ followed by a cyclic sweep of projections
onto the half-spaces via $\mathcal{A}$ as given in \eqref{eq:amsstep}
and \eqref{eq:amssweep}. 

This process was repeated, until the iterate $x^{k}$ became feasible
up to the infeasibility tolerance of $\varepsilon=10^{-8}$ and the
relative change $\frac{\Vert x^{k}-x^{k-1}\Vert}{\Vert x^{k-1}\Vert}$
from the previous iterate was negligible, i. e., became smaller than
$\varepsilon$. Algorithm \ref{supalg} is the pseudo-code for this LinSup
process. All runs were initialized at the all-zeros vector $x^{0}=0$.
The parameters $\varepsilon$, $\alpha$, $\eta_{0}$, $\lambda$, $\tau_{\text{restart}}$
can be adjusted for individual preferences or a particular problem.

\RestyleAlgo{ruled} \LinesNumbered \SetKwInOut{KwParam}{Parameters}
\begin{algorithm}
	\caption{The LinSup superiorization algorithm}
	\label{supalg} 
	\KwData{$A$, $b$, $c$, $\ell$, $u$, $x^{0}$} \KwOut{Superiorized
		vector $y$} \KwParam{$\varepsilon=10^{-8}$, $\alpha=0.99$, $\eta_{0}=10$,
		$\lambda=1$, $\tau_{\text{restart}}=20$} $k\gets0$\; $\tau\gets0$\;
	$\rho\gets0$\; $x^{-1}\gets x^{0}+\mathbbm{1}$\; \While{$\max_{i}(\langle a^{i},x^{k}\rangle-b_{i})\geq\varepsilon$
		or $\frac{\Vert x^{k}-x^{k-1}\Vert}{\Vert x^{k-1}\Vert}\geq\varepsilon$}{
		$x^{k+1}\gets x^{k}-\eta_{k}\frac{c}{\Vert c\Vert}$\; $x^{k+1}\gets\mathcal{A}(x^{k+1})$\;
		$x^{k+1}\gets\max(x^{k+1},\ell)$\; $x^{k+1}\gets\min(x^{k+1},u)$\;
		
		$\tau\gets\tau+1$\; \eIf{$\tau<\tau_{\text{restart}}$}{ $\eta_{k+1}\gets\eta_{k}\alpha$\;
		} { $\tau\gets0$\; $\rho\gets\rho+1$\; $\eta_{k+1}\gets\eta_{0}\alpha^{\rho}$\;
		}
		
		$k\gets k+1$\; } $y\gets x^{k}$\; \Return $y$ 
\end{algorithm}

\subsection{The LP solvers}

\label{optalgs} We compared the actions of LinSup with those of the
following LP solvers: 
\begin{enumerate}
	\item scipy.simplex 
	\item scipy.revised simplex 
	\item scipy.interior-point 
	\item gurobi.primal simplex. 
\end{enumerate}
Scipy is a library for scientific computing in the programming language
Python. It is freely available at \url{www.scipy.org}. Its optimization
suite \textit{scipy.optimize} contains multiple implementations of
common optimization algorithms like the SQP method, the dogleg method
or the conjugate gradient method. Since we are working with the data
of an LP, we employed specialized LP methods only, in particular the
simplex and revised simplex method as well as the interior-point method
of scipy's \textit{linprog} function.

Gurobi is a commercial state-of-the-art solver for linear and nonlinear,
continuous and (mixed) integer problems. It is available at \url{www.gurobi.com}
and offers a wide range of customization, including multiple algorithms,
globalization strategies, preconditioning, etc. For our experiments
we turned off Gurobi's automated choice of algorithm and instead forced
it to use its implementation of the primal simplex to ensure consistency
of the output data.

\subsection{Setting up the experiments}

\label{expsetup} To compare the numerical performance of the different
algorithms with each other in a fair manner we performed the following
steps: First we ran the LinSup algorithm with an infeasibility threshold
of $\varepsilon=10^{-8}$ as a stopping rule. We recorded the runtimes, objective function values
and the maximal constraint violations at each iteration $k$ as well as the total runtime
$T_{\text{sup}}$, i.~e. the time that the superiorization algorithm required until its usual termination criterion. Note that this value will be different for each individual experiment.

Then we turned toward the LP solvers: We imposed on each solver the
same infeasibility threshold of $\varepsilon=10^{-8}$ as a stopping
rule. Additionally we required the solver to stop, if it did not stop
earlier, after the time $T_{\text{sup}}$ was exceeded. In this way,
if the solver would not terminate according to the infeasibility threshold
we stopped it after the time $T_{\text{sup}}$ was reached and then
extracted the current iterate. 

We then compared objective function values and maximal constraint
violation of the LP solver with the corresponding values from the
LinSup algorithm. Similarly, we also documented the runtimes, objective
function values and maximal constraint violations for each iteration
of the LP solvers and plotted the runtimes against objective function
values and maximal constraint violations.

This experimental set up does not aim at reaching LP optimality because
LinSup is not an LP solver. The aim is to reach a certain level of
infeasibility (i. e., getting close to feasibility as we wish) and,
at such a level, compare LinSup with the LP solvers with respect to
the runtimes needed to reach this level of infeasibility and by how
they performed in reducing (not necessarily reaching optimality) the
objective function value.

\section{Experimental Results of the Comparative Investigation}\label{results} 

We split the presentation of our results into two parts, first considering the effects of increasing problem sizes and second considering the effects of increasing condition numbers for fixed sizes.

\begin{figure}
	\begin{subfigure}{0.99\textwidth} \includegraphics[width=0.5\textwidth]{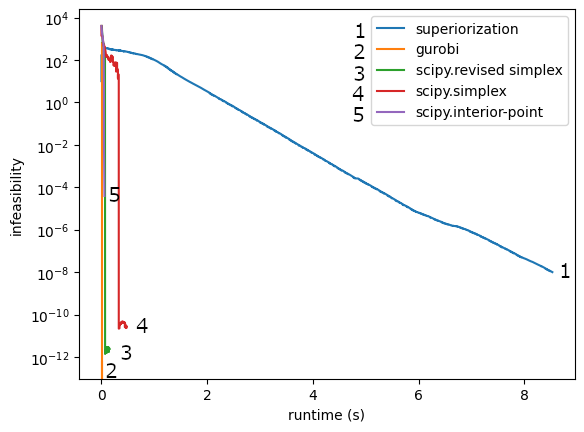}
		\includegraphics[width=0.5\textwidth]{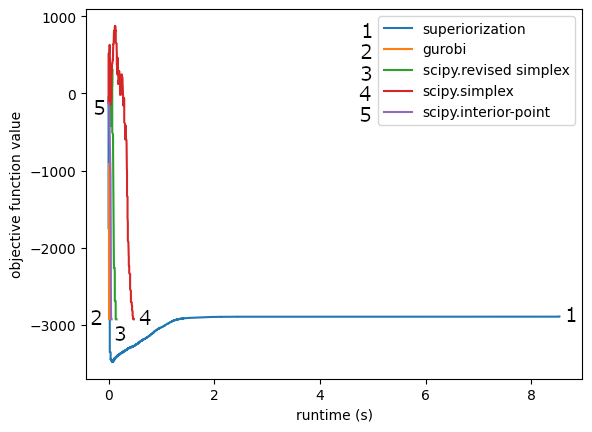}
		\caption{$80\times100$}
		\label{fig:Dim80} \end{subfigure} \begin{subfigure}{0.99\textwidth}
		\includegraphics[width=0.5\textwidth]{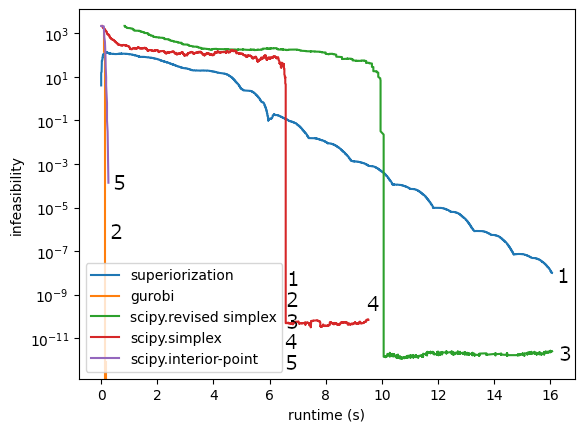}
		\includegraphics[width=0.5\textwidth]{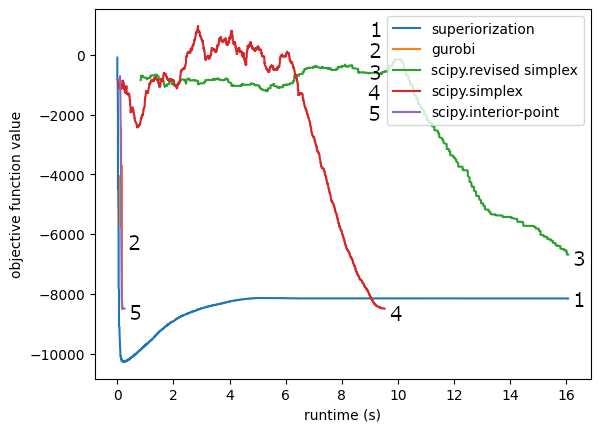}
		\caption{$200\times250$}
	\end{subfigure} 
	\begin{subfigure}{0.99\textwidth} \includegraphics[width=0.5\textwidth]{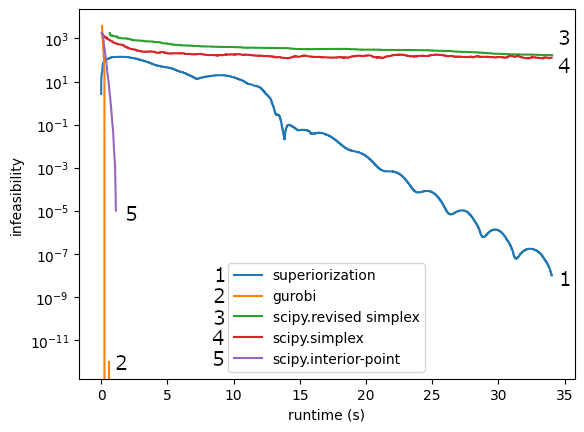}
		\includegraphics[width=0.5\textwidth]{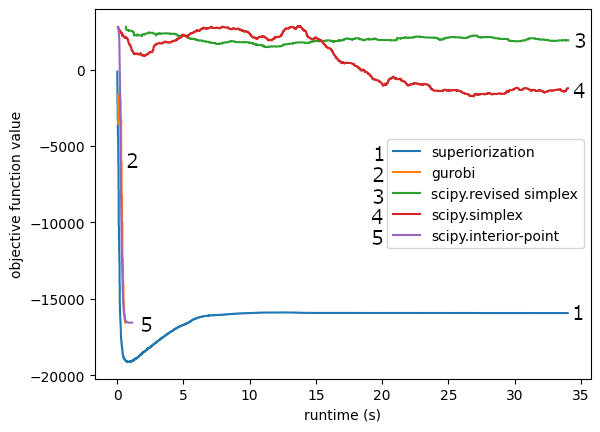}
		\caption{$400\times500$. In this dimension the superiorization algorithm starts
			terminating before scipy's simplex algorithms. We can see that at
			the time of termination, they have not even attained feasibility.}
		\label{fig:Dim400} \end{subfigure} 
	\caption{Dimensions $80\times100$ -- $400\times500$ at fixed $\kappa=1000$.}
\end{figure}

\begin{figure}
	\begin{subfigure}{0.99\textwidth} \includegraphics[width=0.5\textwidth]{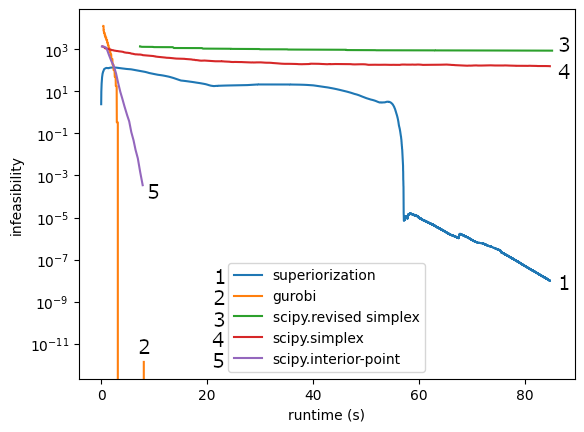}
		\includegraphics[width=0.5\textwidth]{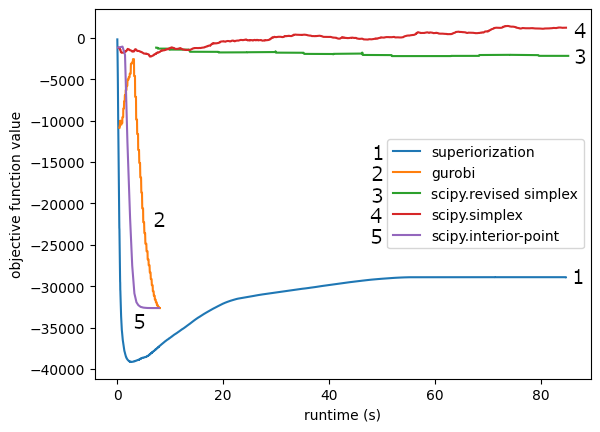}
		\caption{$800\times1000$}
		\label{fig:Dim800} \end{subfigure}
	
	\begin{subfigure}{0.99\textwidth} \includegraphics[width=0.5\textwidth]{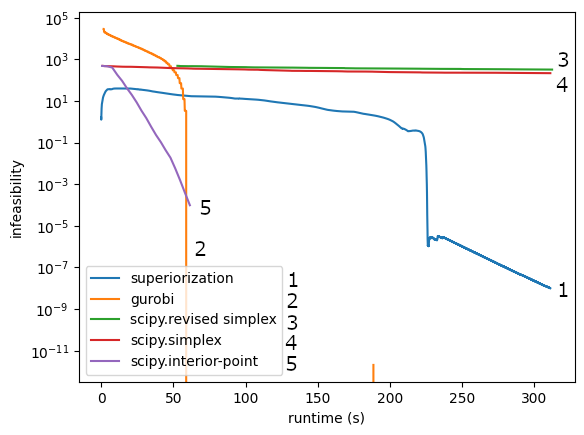}
		\includegraphics[width=0.5\textwidth]{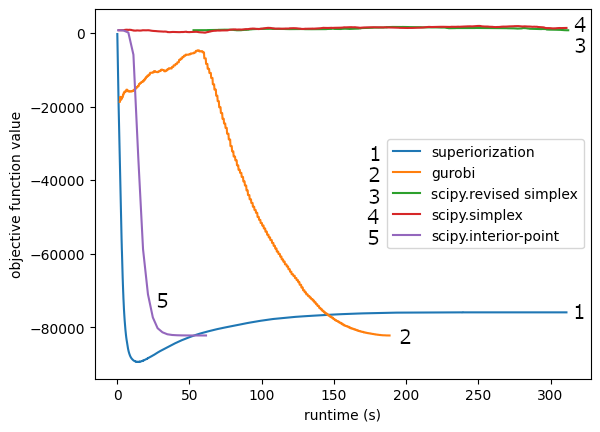}
		\caption{$2000\times2500$}
		\label{fig:Dim2000} \end{subfigure}
	
	\begin{subfigure}{0.99\textwidth} \includegraphics[width=0.5\textwidth]{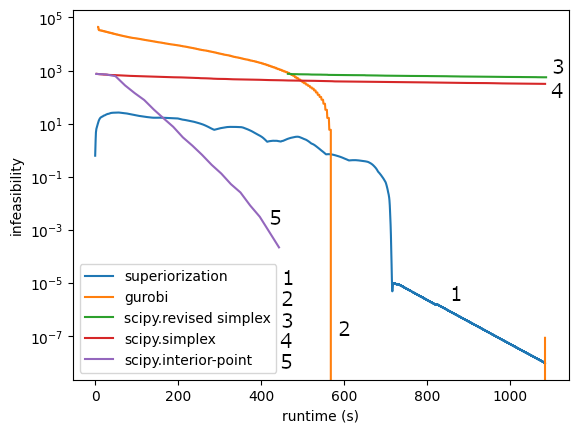}
		\includegraphics[width=0.5\textwidth]{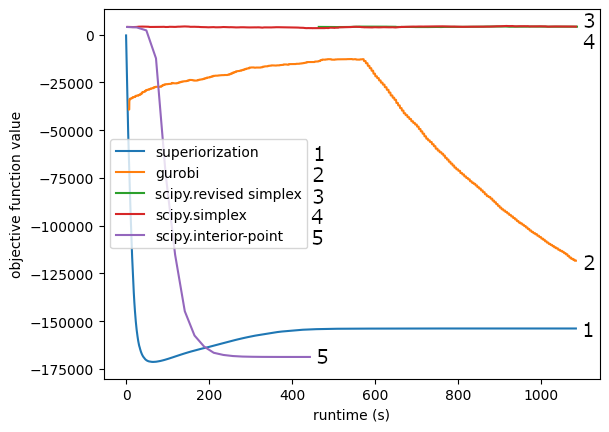}
		\caption{$4000\times5000$. In this dimension we see that the superiorization
			algorithm terminates even before Gurobi's simplex does.}
		\label{fig:Dim4000} \end{subfigure} 
	\caption{Dimensions $800\times1000$ -- $4000\times5000$ at fixed $\kappa=1000$.} 
\end{figure}

\subsection{Sensitivity to increasing sizes when the condition number is fixed}

Figures \ref{fig:Dim80}--\ref{fig:Dim4000} show the behavior of
the different algorithms on problems with increasing dimension for
fixed condition number $\kappa=1000$. On the left we see the maximum
violation of the constraints $\max_{i}(\langle a^{i},x\rangle-b_{i})$
plotted against the runtime. On the right we have the corresponding
objective function values $\langle c,x\rangle$ plotted against the
runtime. 

The trend is clear: While in problems with smaller dimensions the
simplex algorithms are considerably faster, with increasing problem
dimensions the simplex implementations take much longer to terminate,
up to the point where the superiorization method terminates quicker.
We also see a clear difference in the algorithms styles: Simplex aims
to become feasible first and then starts to improve the objective
function. Superiorization, on the other hand, reduces objective function
values in its initial iterations because then the step-sizes of the
perturbations are still large and only as iterations proceed the effect
of feasibility-seeking becomes stronger.
Further,
Gurobi outperformed both of the other simplex implementations, which
is to be expected by a commercial product. 

Nevertheless, we can see
that Gurobi, too, needed more and more time until termination with
increasing dimension and eventually the superiorization algorithm
terminated quicker.
Summarizing, our first finding is, that superiorization reaches its termination criterion
considerably faster than the simplex implementations in higher dimensions. This shows that superiorization can produce
better results than the simplex implementations, when both are terminated
at a certain time, that is before the usual termination criterion
is met. 
Furthermore, based on our computational
experiments, we see that the superiorization methodology is quite robust with
respect to increasing dimensions. 

In particular, the runtime
remains relatively small when compared to the scipy implementations
of the simplex algorithm.


\subsection{Immunity to increased condition numbers}

\begin{figure}
	\begin{subfigure}{0.99\textwidth} \includegraphics[width=0.5\textwidth]{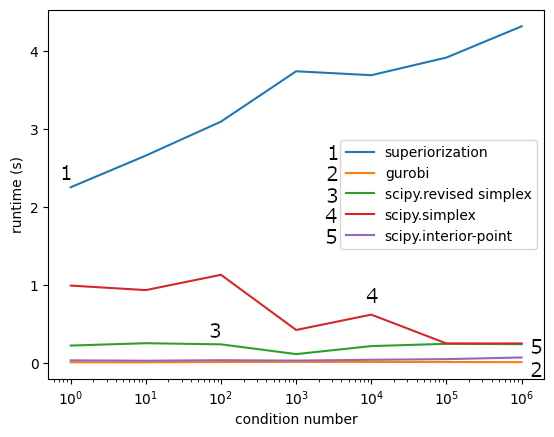}
		\includegraphics[width=0.5\textwidth]{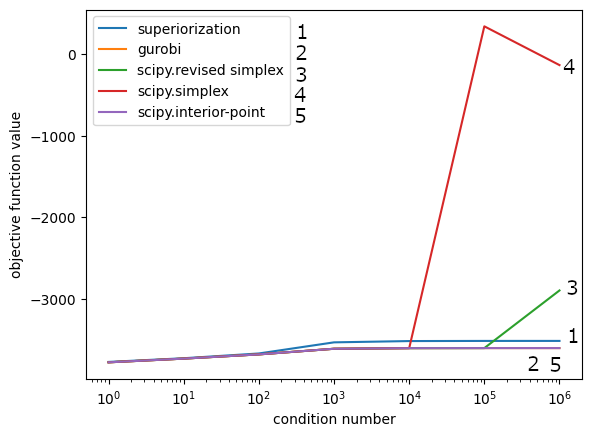}
		\caption{$80\times100$}
		\label{fig:accDim80} \end{subfigure}  
	\begin{subfigure}{0.99\textwidth} \includegraphics[width=0.5\textwidth]{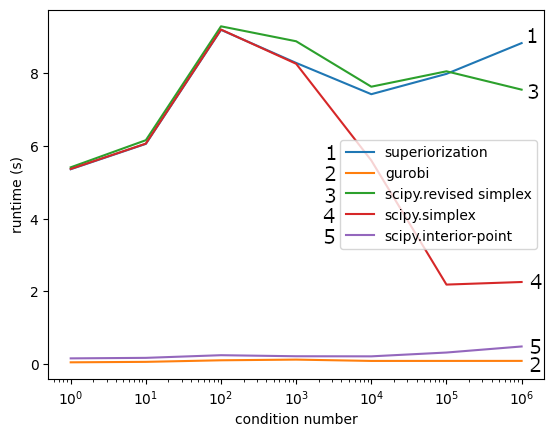}
		\includegraphics[width=0.5\textwidth]{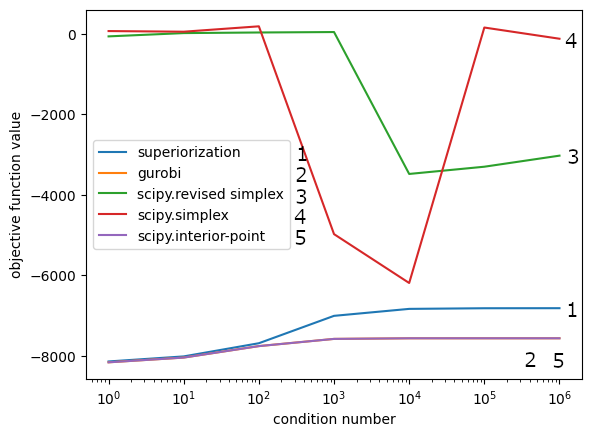}
		\caption{$200\times250$}
		\label{fig:accDim200} \end{subfigure}
	\begin{subfigure}{0.99\textwidth} \includegraphics[width=0.5\textwidth]{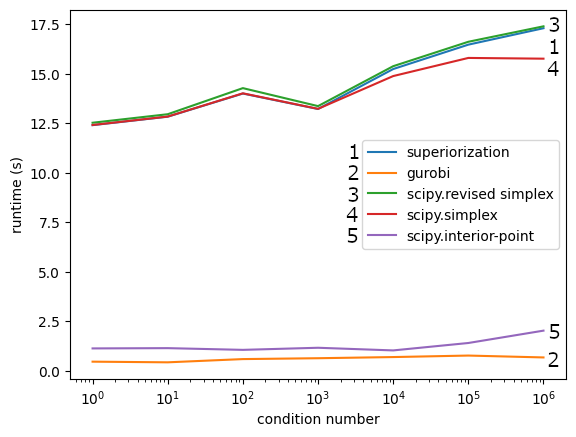}
		\includegraphics[width=0.5\textwidth]{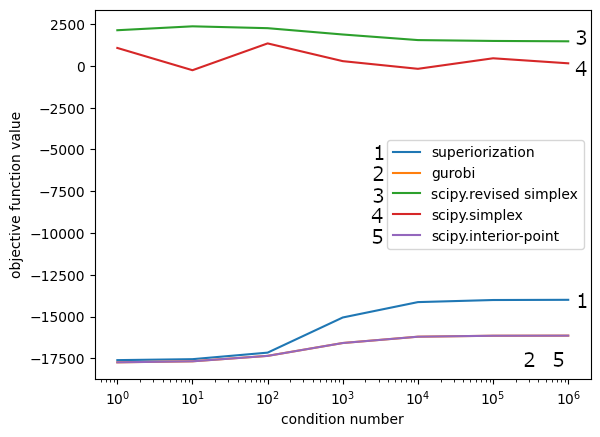}
		\caption{$400\times500$}
		\label{fig:accDim400} \end{subfigure}
	
	\caption{Behavior of the algorithms at fixed problem sizes $80\times 100$ -- $400\times 500$ and increasing condition numbers.}
\end{figure}

\begin{figure}
	\begin{subfigure}{0.99\textwidth}
		\includegraphics[width=0.5\textwidth]{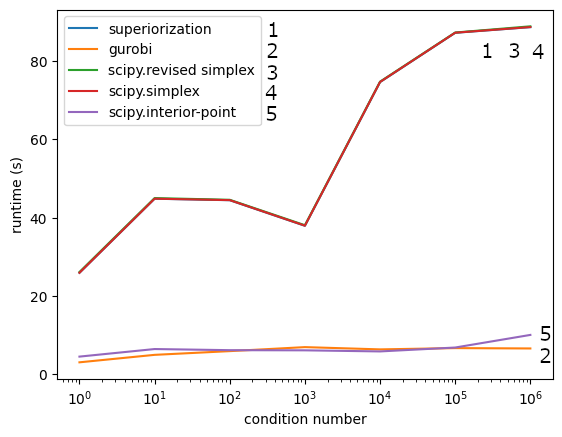} \includegraphics[width=0.5\textwidth]{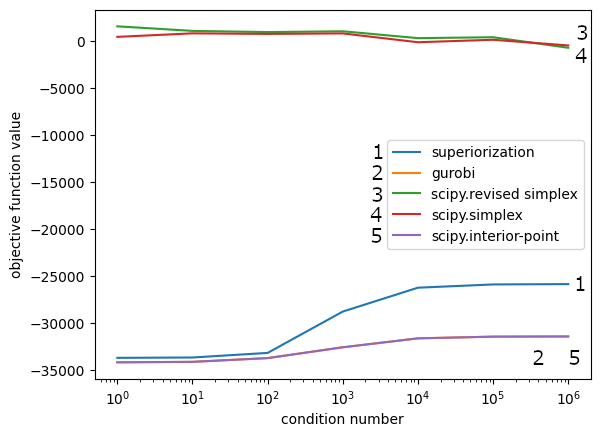}
		\caption{$800\times1000$}
		\label{fig:accDim800} 
	\end{subfigure}
	\begin{subfigure}{0.99\textwidth}
		\includegraphics[width=0.5\textwidth]{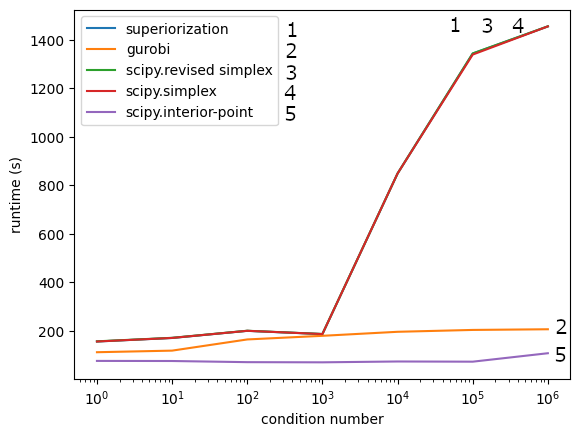} \includegraphics[width=0.5\textwidth]{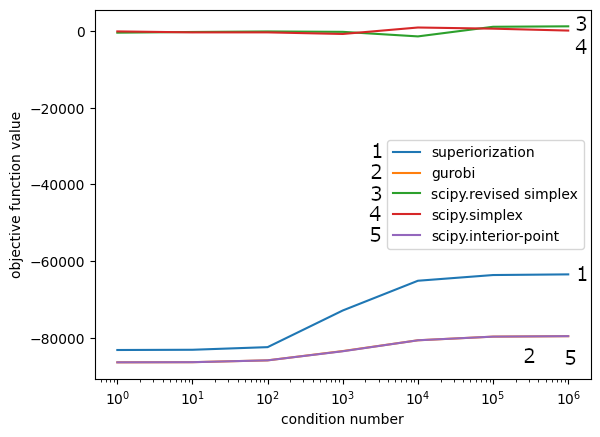}
		\caption{$2000\times2500$}
		\label{fig:accDim2000} 
	\end{subfigure} 
	\begin{subfigure}{0.99\textwidth}
		\includegraphics[width=0.5\textwidth]{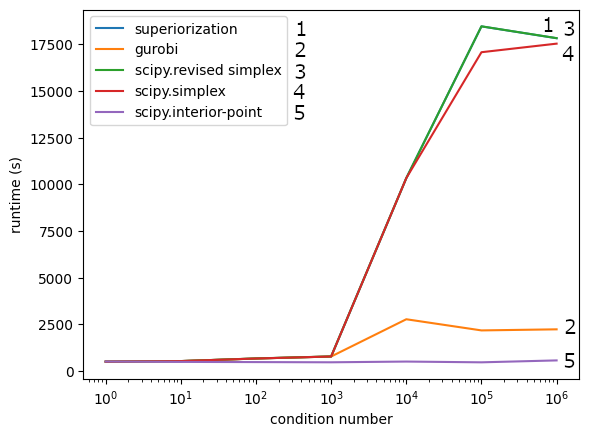} \includegraphics[width=0.5\textwidth]{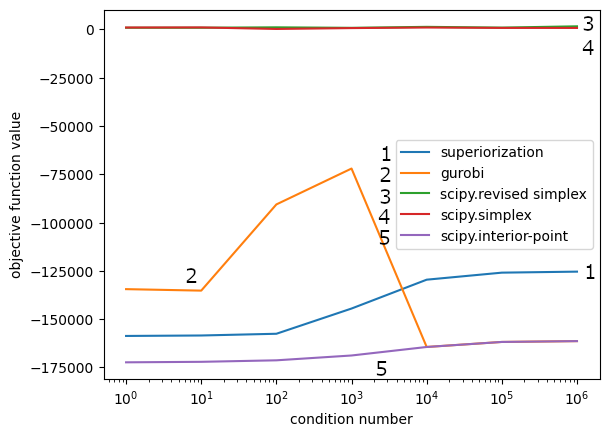}
		\caption{$4000\times5000$}
		\label{fig:accDim4000} 
	\end{subfigure}
	\caption{Behavior of the algorithms at fixed problem sizes $800\times 1000$ -- $4000\times 5000$ and increasing condition numbers.}
\end{figure}

We move on to the evaluation of the experiments with increasing condition numbers. Figures \ref{fig:accDim80} -- \ref{fig:accDim4000} demonstrate well
the robustness with respect to increasing condition numbers. On the left we see the averaged runtimes of $5$ problem instances for each algorithm, plotted against the condition numbers, on the right we see the averaged objective function values at termination, plotted against the condition numbers. 

Notice the severely sub-optimal
objective function value of scipy's simplex for high condition numbers. This explains the low runtime as the algorithm realizes that these problems are hard and quickly ``surrenders".

The trend continues in Figures \ref{fig:accDim800} -- \ref{fig:accDim4000} and the revised simplex too starts to reach its limits. Superiorization on the other hand proves to be quite stable in terms of its runtime with regards to increasing condition numbers, at times surpassing the Gurobi implementation in both runtime and objective function value. While the total runtime \textit{does} increase significantly with the condition number, it remains within reasonable time frames, considering the increased difficulty of the problems that accompanies higher condition numbers.

\begin{figure}
	\begin{subfigure}{0.495\textwidth}
		\includegraphics[width=\textwidth]{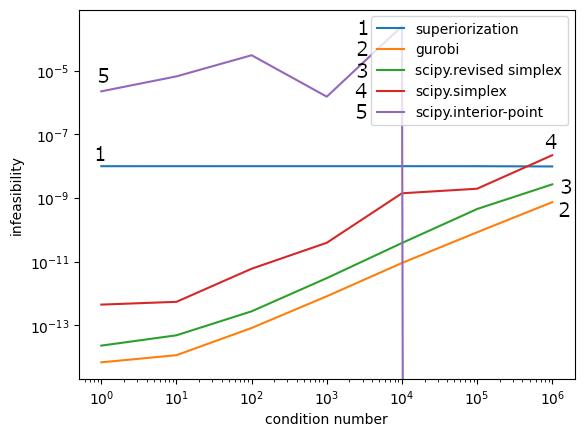}
		\caption{$80\times100$}
		\label{fig:accInfDim80} 
	\end{subfigure}
	\begin{subfigure}{0.495\textwidth}
		\includegraphics[width=\textwidth]{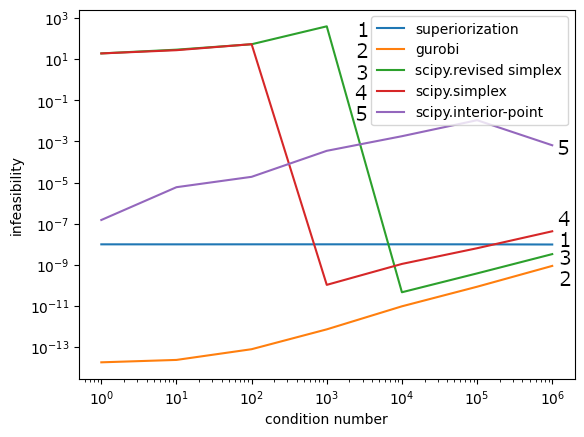} 
		\caption{$200\times250$}
		\label{fig:accInfDim200} 
	\end{subfigure} 
	\begin{subfigure}{0.495\textwidth}
		\includegraphics[width=\textwidth]{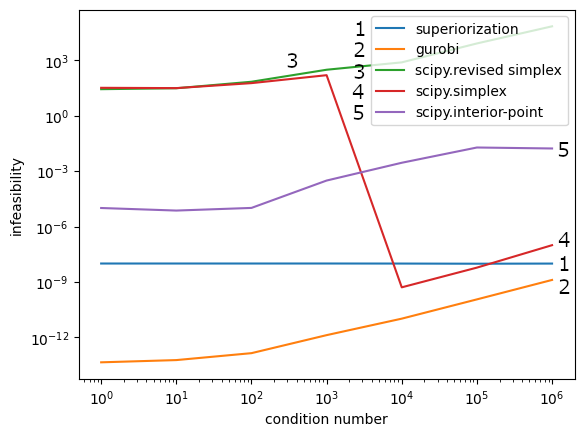} 
		\caption{$400\times500$}
		\label{fig:accInfDim400} 
	\end{subfigure}
	\begin{subfigure}{0.495\textwidth}
		\includegraphics[width=\textwidth]{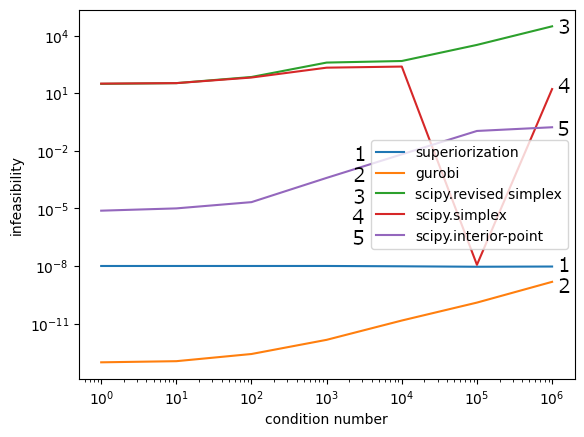} 
		\caption{$800\times1000$}
		\label{fig:accInfDim800} 
	\end{subfigure}
	\begin{subfigure}{0.495\textwidth}
		\includegraphics[width=\textwidth]{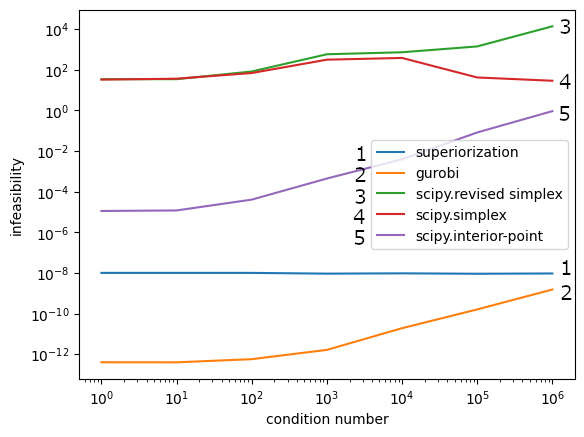} 
		\caption{$2000\times2500$}
		\label{fig:accInfDim2000} 
	\end{subfigure}
	\begin{subfigure}{0.495\textwidth}
		\includegraphics[width=\textwidth]{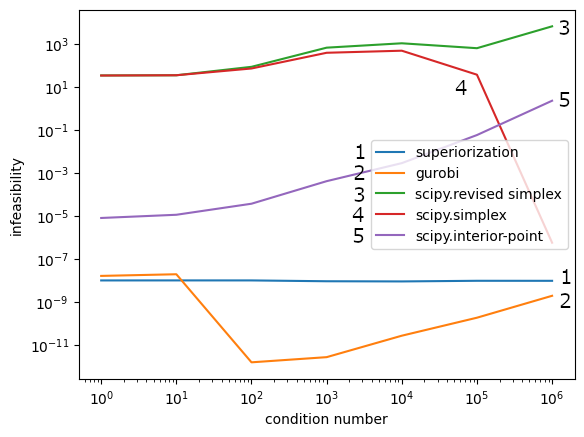} 
		\caption{$4000\times5000$}
		\label{fig:accInfDim4000} 
	\end{subfigure}
	\caption{Average infeasibilities at termination of the problem instances. }
\end{figure}

\begin{figure}
	\begin{subfigure}{0.495\textwidth}
		\includegraphics[width=\textwidth]{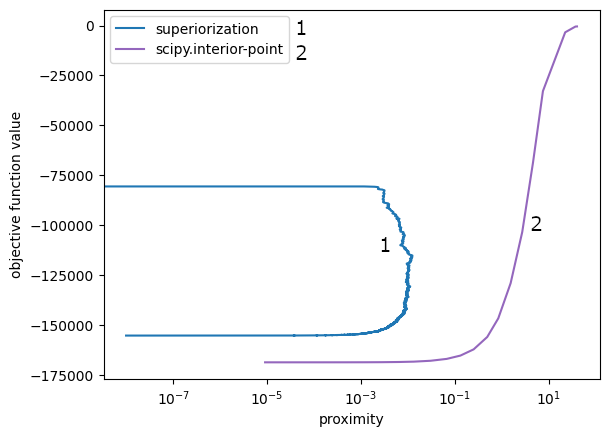}
		\caption{$\kappa=1$}
		\label{fig:proxTarg1} 
	\end{subfigure}
	\begin{subfigure}{0.495\textwidth}
		\includegraphics[width=\textwidth]{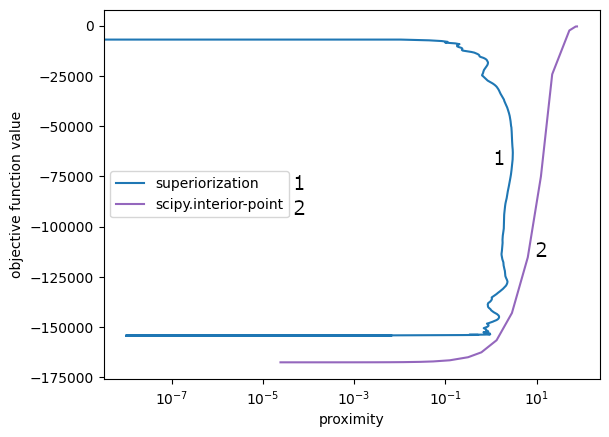} 
		\caption{$\kappa=10^2$}
		\label{fig:proxTarg100} 
	\end{subfigure} 
	\begin{subfigure}{0.495\textwidth}
		\includegraphics[width=\textwidth]{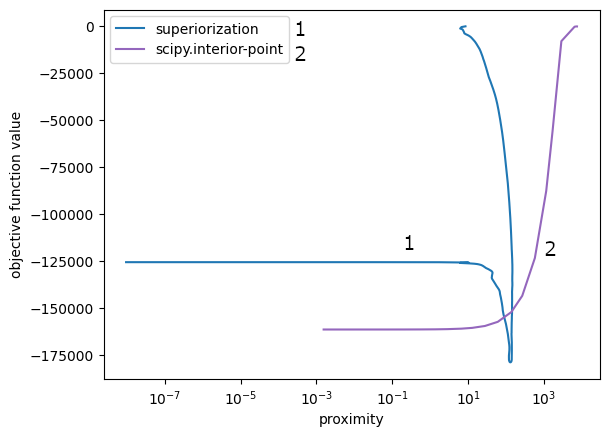} 
		\caption{$\kappa=10^4$}
		\label{fig:proxTarg10000} 
	\end{subfigure}
	\begin{subfigure}{0.495\textwidth}
		\includegraphics[width=\textwidth]{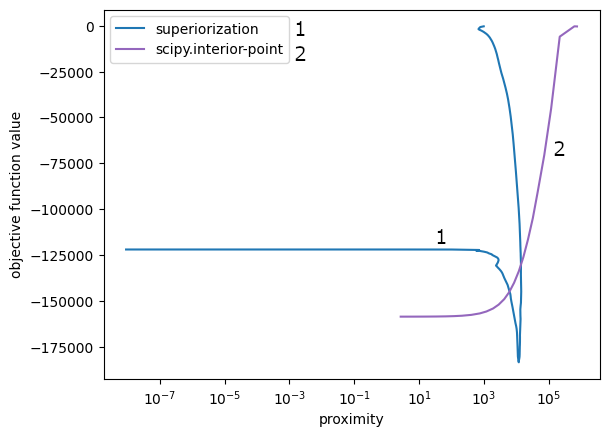} 
		\caption{$\kappa=10^6$}
		\label{fig:proxTarg1000000} 
	\end{subfigure}
	\caption{Proximity-Target-Curves of superiorization and of the interior-point method for multiple condition numbers. Notice once again the severe infeasibility (i. e., proximity) of some of the final iterates of the interior-point method.}
\end{figure}

Next, we take a look at the infeasibility at termination (Figures \ref{fig:accInfDim80}-\ref{fig:accInfDim4000}). Notice how with increasing condition number the general trend is towards higher infeasibility, except for superiorization, which remains constant throughout. Like before, scipy.simplex and scipy.revised simplex turn out to be unreliable in moderate dimension already. Additionally, we see that the interior point method has trouble attaining the desired level of feasibility. 

We can see values of infeasibility of magnitudes up to $10^0-10^1$, significantly increasing with the condition number. In other words, while the interior point method seems to be able to quickly find the correct objective function value, it has issues doing so in a feasible manner, making it unreliable for practical applications with high condition numbers. 

On the other hand, superiorization consistently reaches precisely the required level of feasibility, regardless of the condition number. This becomes particularly evident when looking at the proximity-target-curves, i.~e. the plot of the level of infeasibility against the objective function value for all iterates of the respective algorithms. For further information on proximity-target curves for evaluating performance of iterative algorithms see \cite[Section 4]{censorprox}. Figures \ref{fig:proxTarg1}-\ref{fig:proxTarg1000000} showcase the behavior of the interior point method as well as of superiorization.

Hence, our second finding is, that the superiorization methodology is able to easily handle large condition numbers, especially when compared to the scipy implementations. We conjecture that this is due to the fact,
that our superiorization implementation never considers the full problem
at once, but performs individual projections onto the half-spaces
instead. This comes at the cost of sacrificing feasibility during
early iterations. 

Furthermore, we conjecture, that the bounded perturbation
resilience of the basic algorithm may play a role in absorbing errors
that occur during computation, which would normally be amplified by
the condition number. 


This would mean that superiorization, in general,
may be less affected by high condition numbers, than other current
algorithms.

We observe that in our experiments, in terms of runtime and objective function value, at first glance the interior point method seemed to outperformed
the other algorithms. This can
be explained in the following way: 

The condition number of a matrix
can be interpreted as a measure of how linearly dependent its rows
or columns are. A well-conditioned matrix (i.~e., $\kappa=1$) only
has a single singular value and will be semi-orthogonal, whereas an
ill-conditioned matrix ($\kappa$ ``large'') will have ``almost''
linearly dependent entries. Consequently, the half-spaces $\{x\in\mathbb{R}^{n}\colon\langle a^{i},x\rangle\leq b_{i}\}$
will be almost parallel and the resulting polyhedron will consist
of many facets and vertices. 

A basic simplex implementation, which
moves from vertex to vertex, will consequently face a long runtime.
The interior point method, on the other hand, is not dependent on
the vertices. It will take its path through the interior of the polyhedron
regardless of its boundary.\\

However, as stated above, scipy's interior point method seems to have trouble finding a feasible point at all for larger condition numbers, showcasing the strength of superiorization.

Combining all of the above observations, superiorization clearly offers a valuable contribution to large-scale optimization problems of high condition numbers. For large-enough problem instances it is faster than scipy's simplex methods and more reliable than scipy's interior point method in terms of feasibility.

\section{Conclusions}

\label{conclusion} In this paper we evaluate experimentally the
superiorization method (SM) and constrained optimization algorithms on a set of exemplary
linear problems with varying condition numbers and varying problem sizes with
the aim of investigating and comparing their immunity to increasing
condition numbers. These results expand on the results in \cite{linsup} and further show the strengths and weaknesses of linear superiorization and classical LP optimization algorithms. Our experimental results are promising for the
observed problem sizes, and we are confident that the trend that we
observed will continue for larger problem instances.

The superiorization method and constrained optimization algorithms use the same
input data which consists of a family of constraints obtained from
the modeling process along with a user-chosen objective function.
But the two approaches aim at different end-points of their iterative
processes. The easy accessibility of the superiorization methodology
allows for quick implementations with the advance knowledge that the
aim is not to reach a constrained optimum. At the same time the SM
may compute its solutions at a lower runtime (in the case of simplex) and give better feasibility than the interior-point method.

As superiorization is a relatively new concept, we expect that, with
further tuning of its parameters it will continue to find a place
as a computational model and tool in situations in which users do
not wish to invest efforts in seeking a constrained optimal point
but rather wish to find a feasible point which is “superior” in the
sense of having a smaller or equal objective function value than that
of a feasible point reached by the same feasibility-seeking algorithm.

Another interesting point is that, as is well-known, interior point
methods reach their performance limits for ill-conditioned nonlinear
problems (e.~g., in intensity-modulated radiotherapy). It would be
interesting to compare the superiorization methodology in a nonlinear
setting, with a different basic feasibility-seeking algorithm, to
the interior point method to see if superiorization can contribute
to solving these problems faster.

\vskip 6mm
\noindent{\bf Acknowledgments}

\noindent The authors thank the anonymous referee for valuable comments and suggestions, which helped improve the presentation and clarity of the paper.\\

\noindent  The work of Y.C. is supported by U.S. National
Institutes of Health Grant Number R01CA266467 and by the Cooperation
Program in Cancer Research of the German Cancer Research Center (DKFZ)
and Israel's Ministry of Innovation, Science and Technology (MOST).

\bibliography{bibliography}
\bibliographystyle{myunsrt}

	

\end{document}